\documentclass[trsc,nonblindrev]{informs3} 

\OneAndAHalfSpacedXI 

\renewcommand{\theARTICLETOP}{}
\renewcommand{\theRRHFirstLine}{}
\renewcommand{\theRRHSecondLine}{}
\renewcommand{\theLRHFirstLine}{}
\renewcommand{\theLRHSecondLine}{}
\usepackage{soul}
\usepackage{natbib}
 \bibpunct[, ]{(}{)}{,}{a}{}{,}%
 \def\bibfont{\small}%
 \usepackage{mathtools}
\usepackage{mathpazo}
\usepackage{amsmath}
\usepackage{mathrsfs}
\usepackage{amssymb}
\usepackage{multirow}
\usepackage{longtable}
\usepackage{graphicx}
\usepackage{enumerate}
\usepackage[titletoc]{appendix}
\usepackage{graphicx} 
\usepackage{epstopdf}
\usepackage{subfigure}
\usepackage{algorithm,algorithmic,longtable}
\usepackage[bookmarks=true,linkcolor=red,citecolor=blue,urlcolor=green,colorlinks=true, breaklinks]{hyperref}
\usepackage{indentfirst}

\usepackage{booktabs}
\usepackage{array, caption, threeparttable,arydshln}
\usepackage[font=footnotesize,labelfont=bf,labelsep=none]{caption}
\usepackage{caption}
\DeclareCaptionLabelSeparator{twospace}{\ ~\ ~}   
\newtheorem{myThe}{Theorem}

\newtheorem{myCor}{Corollary}

\usepackage{empheq}
\usepackage{cases}

\allowdisplaybreaks[1]
\usepackage{longtable}
\usepackage{multicol} 
\renewcommand{\algorithmicrequire}{\textbf{Input:}}
\renewcommand{\algorithmicensure}{\textbf{Output:}}

\makeatletter
\newcommand{\mathleft}{\@fleqntrue\@mathmargin\parindent}
\newcommand{\mathcenter}{\@fleqnfalse}
\makeatother
\usepackage{etoolbox}
\newcommand{\zerodisplayskips}{%
\setlength{\abovedisplayskip}{7pt}%
\setlength{\belowdisplayskip}{6pt}%
\setlength{\abovedisplayshortskip}{7pt}%
}
\appto{\normalsize}{\zerodisplayskips}
\appto{\small}{\zerodisplayskips}
\appto{\footnotesize}{\zerodisplayskips}

\allowdisplaybreaks \allowdisplaybreaks[1]

\TheoremsNumberedThrough     

\EquationsNumberedThrough    

\MANUSCRIPTNO{} 

\begin{document}

\TITLE{Modular Vehicle Routing Problem: Applications in Logistics}

\ARTICLEAUTHORS{
\AUTHOR{Hang Zhou, Yang Li, Chengyuan Ma, Keke Long, Xiaopeng Li\footnote{Corresponding author. Email: xli2485@wisc.edu}}
\AFF{Department of Civil and Environmental Engineering, University of Wisconsin-Madison}
} 

\ABSTRACT{
    Recent studies and industry advancements indicate that modular vehicles (MVs) have the potential to enhance transportation systems through their ability to dock and split during a trip. Although various applications of MVs have been explored across different domains, their application in logistics remains underexplored. This study examines the use of MVs in cargo delivery to reduce total delivery costs. We model the delivery problem for MVs as a variant of the Vehicle Routing Problem, referred to as the Modular Vehicle Routing Problem (MVRP). In the MVRP, MVs can either serve customers independently or dock with other MVs to form a platoon, thereby reducing the average cost per unit. In this study, we mainly focus on two fundamental types of MVRPs, namely the capacitated MVRP (CMVRP) and the MVRP with time windows (MVRPTW). To address these problems, we first developed mixed-integer linear programming (MILP) models, which can be solved using commercial optimization solvers. Given the NP-hardness of this problem, we also designed a Tabu Search (TS) algorithm with a solution representation based on Gantt charts and a neighborhood structure tailored for the MVRP. Multi-start and shaking strategies were incorporated into the TS algorithm to escape local optima. Additionally, we explored other potential applications in logistics and discussed problem settings for three MVRP variants. Results from numerical experiments indicate that the proposed algorithm successfully identifies nearly all optimal solutions found by the MILP model in small-size benchmark instances, while also demonstrating good convergence speed in large-size benchmark instances. Comparative experiments show that the MVRP approach can reduce costs by approximately 5.6\% compared to traditional delivery methods. Sensitivity analyses reveal that improving the cost-saving capability of MV platooning can enhance overall benefits.
}

\maketitle

\section{Introduction}

The rapid growth of e-commerce, local commerce, and the retail industry necessitates enhanced service capabilities in city-scale logistics \citep{janjevic2019integrating, zhou2023exact}. Traditional city-scale logistics often rely on a mixed fleet of vehicles, where large trucks transport significant volumes of cargo to central distribution hubs, followed by smaller vehicles managing last-mile delivery to localized distribution points within communities \citep{zhou2022two}. However, this conventional model requires logistics providers to maintain and allocate a diverse fleet based on cargo volumes, leading to increased costs related to labor, time, and infrastructure due to cargo transfers between different vehicle types at specific stations.

Fortunately, with the emergence of modular vehicle (MV) technologies, flexible docking, splitting, and even cargo transfer between vehicles en route are now feasible \citep{chen2022continuous, shi2021operations, zhang2020modular}. The defining feature of MVs is their modularity, which allows them to dock and split en route to form platoons of varying sizes. Utilizing MVs for delivery eliminates the need for a diverse vehicle fleet, thereby reducing vehicle costs. Furthermore, MV platoons can adjust their capacity en route, eliminating additional time costs associated with cargo transfers. Unlike traditional vehicle platoons, MVs connect physically to form platoons, significantly reducing fuel consumption \citep{li2022review}. Therefore, the application of MVs in logistics, particularly in last-mile delivery, offers numerous advantages \citep{hannoun2022modular}.

There has been some discussion in academia regarding the application of MVs. Most studies focus on the use of MVs in transit systems, such as shuttle systems \citep{chen2019operational, chen2020operational}, corridor systems \citep{shi2021operations, chen2022continuous}, and taxi systems \citep{fu2023dial}. Recent research demonstrates the significant potential of MVs in transit systems, showing a reduction in total costs by 9.87\%–32.09\% \citep{zhang2024optimising}. A few studies have explored the application of MVs in Emergency Medical Services (EMS) \citep{hannoun2022modular} and joint passenger and freight transport \citep{hatzenbuhler2023modular}, confirming the benefits of MVs in these contexts as well. Moreover, The concept of MVs is already being applied in real-world scenarios. Next Future Transportation Inc. (\url{https://www.next-future-mobility.com/}) is a leader in MV technology development, especially with its modular bus. This bus can flexibly adjust its capacity through docking and splitting capabilities, addressing challenges such as bus bunching due to unpredictable delays \citep{khan2023application}. Although the design and implementation of MVs are still in the early stages, it is crucial for the research community to explore the potential of this innovative technology further. Moreover, the application of MVs in logistics, particularly in last-mile delivery, warrants greater attention.

While numerous studies on MVs have emerged in recent years, most of them focus on designing modular transit systems and overlook the significant benefits that MVs could offer to logistics applications. Integrating MVs into traditional delivery problems presents unique challenges in both modeling and algorithm development. Since MVs can dock together in a platoon with expanded capacity, traditional capacity constraints can not be applied to MVs. The capacity of a platoon is the sum of the capacity of all MVs within it. Consequently, the capacity of a platoon changes dynamically as MVs dock or split, necessitating constant tracking of the spatial and temporal status of the MVs.

To address this research gap, our study proposes a novel variant of the Vehicle Routing Problem (VRP), referred to as the Modular Vehicle Routing Problem (MVRP). To the best of our knowledge, this is the first work to introduce the MVRP in logistics. In the MVRP, MV platoons depart from a distribution center to serve a set of customers. Each MV can operate independently or dock with others at customer locations to form a platoon, which can subsequently split into smaller platoons if needed. We primarily consider two fundamental types of the MVRP, namely the capacitated MVRP (CMVRP) and the MVRP with time windows (MVRPTW), which are variants of the well-known capacitated VRP and VRP with time windows, respectively. Our research analyzes these problems in three aspects. First, theoretical analyses prove that the optimal solutions of the MVRP have upper and lower bounds relative to the optimal solutions of the VRP. Second, to solve the CMVRP and MVRPTW, we develop Mixed Integer Linear Programming (MILP) models and a tailored Tabu Search (TS) algorithm, which can optimally solve small-size instances and efficiently handle large-size instances, respectively. Finally, we introduce several potential applications in logistics and propose MVRP problem settings for future research. Our work contributes to the literature in the following ways:
\begin{enumerate}
   \item We introduce a new variant of the VRP that incorporates MV technology, referred to as the MVRP, where delivery vehicles can dock into platoons and split in a trip. The theoretical lower bound of the MVRP is compared with that of the traditional VRP.
   \item We develop three-index MILP models and a TS algorithm to solve the CMVRP and MVRPTW. A tailored TS algorithm is built upon a specialized data structure, the Gantt chart, which is widely used in scheduling algorithms. Several components of the TS algorithm, including a customized neighborhood structure, an acceleration strategy for verifying time window constraints, as well as multi-start and shaking strategies are specifically designed for the MVRP. Numerical experiments demonstrate the effectiveness of the proposed TS algorithm.
   \item The potential logistics application scenarios and problem settings for three MVRP variants are discussed, which provides directions for future research.
\end{enumerate}

The remainder of the paper is as follows. Section \ref{sec:literature} reviews related studies. Section \ref{sec:problem} defines the problem and constructs the MILP model. Section \ref{sec:algorithm} presents the TS algorithm. Section \ref{sec:variants} discusses several variants of the MVRP. Section \ref{sec:experiments} describes the instance sets and provides the numerical results. Section \ref{sec:conclusion} concludes this paper and suggests future research directions.

\section{Literature Review}
\label{sec:literature}

Given that MV is an emerging technology, there is limited literature specifically focusing on its application in logistics and routing problems. Therefore, in this section, we primarily review the application of MVs and the related VRP variants in the literature. 

\subsection{Applications of MVs}

A few studies have investigated the application of MVs, particularly in transit systems. \cite{chen2019operational} and \cite{chen2020operational} proposed both discrete and continuous methods for the joint design of dispatch headway and capacity in a one-to-one shuttle system using modular autonomous vehicles (MAVs). \cite{chen2021designing, chen2022continuous} and \cite{shi2021operations} addressed the operational design problem for urban mass transportation corridor systems with MAVs. However, their research primarily focuses on the operational design of fixed-route systems. \cite{tian2022planning} explored the optimal planning of public transit services with MVs, allowing for station-wise docking and splitting. Their study optimizes the location and capacity of stations. Recently, \cite{zhang2024optimising} optimized the MAV transit service employing docking-splitting operations plus the skip-stop strategy. Some articles have focused on the trajectory optimization of MAVs during en route docking and splitting operations, such as \cite{li2022trajectory} and \cite{li2023trajectory}. These studies provide valuable insights for the broader application of MVs in various scenarios.

Although most early studies in MVs consider the transit system with fixed routes, several studies are modeling and solving the routing problems.
\cite{gong2021transfer} investigated a customized modular bus system where the number of modular units in each bus can vary. They proposed a network design problem that jointly optimizes transfer-based customized bus routes and passenger-route assignments for each departure time, presenting a variant of the VRP. The authors developed a MILP model and a Particle Swarm Optimization (PSO) heuristic algorithm. Their experiments demonstrate that the PSO algorithm achieves an average gap of 5.26\% compared to the exact solution.
\cite{hannoun2022modular} introduced a smart EMS system that leverages the docking and splitting operations of MVs for rapid patient and medical personnel transfer. The authors designed a MILP model to solve this problem.
\cite{fu2023dial} introduced a Modular Dial-a-Ride Problem (MDARP). This system enables en-route passenger transfers before the splitting of MVs. The authors considered the soft deadline constraints and the extension of hard time windows. However, the study limits its generalizability by only considering two vehicles per platoon. To effectively address this problem, they propose a heuristic algorithm featuring a Steiner-tree-inspired neighborhood search and an improvement heuristic.
\cite{hatzenbuhler2023modular} explored using MVs for joint passenger and freight transport. The authors introduce a variation of the Pickup and Delivery Problem (PDP), which is considered the MV platoon with a flexible number of units but is predetermined at the depot. They proposed an Adaptive Large Neighborhood Search (ALNS) algorithm whose optimal gap is less than 0.01\% in 30 smaller instances with fewer than 17 nodes. However, as splitting operations are not considered, this problem closely resembles a heterogeneous VRP. 


\subsection{Vehicle Platooning Problem}

A related research problem to the MVRPTW is the Vehicle Platooning Problem (VPP), also referred to in the literature as the Platooning Routing Problem (PRP) or the Truck Platooning Problem (TPP). In the VPP, a fleet of vehicles is given a set of starting points and destinations. These vehicles can form platoons with other vehicles along their routes to save fuel costs. The VPP aims to find the optimal platoon routing for these vehicles on a graph. A comprehensive review of the VPP can be found in \cite{bhoopalam2018planning}.
Starting with \cite{larsson2015vehicle}, many studies have investigated the VPP. \cite{larsson2015vehicle} first formally defined the VPP as a variant of the Vehicle Routing Problem (VRP) that minimizes fuel consumption and demonstrated that this problem is NP-hard. The authors then developed a Mixed Integer Linear Programming (MILP) formulation for the VPP and presented two constructive heuristics and one local search algorithm for large-scale instances.
\cite{luo2018coordinated} proposed a coordinated platooning MILP model that integrates speed selection and platoon formation/dissolution into the problem formulation.
\cite{luo2022repeated} built on the idea that trucks cannot deviate significantly from their shortest paths to form platoons. They proposed an iterative "route first, schedule second" approach to plan platoons and developed a set of valid inequalities for both steps to reduce the problem size.
\cite{bhoopalam2023platoon} focused on the scenario where the maximum platoon size is two. They provided a polynomial algorithm for this case and, based on this, designed two fast heuristics for the multi-truck platooning problem. They also conducted numerical tests on a Dutch highway network consisting of 20 cities and 45 road sections; the results indicate that two-truck platoons can capture most of the potential platooning savings.
\cite{zhao2024improved} considered the constraint that trucks are not allowed to wait during their trip and proposed a heuristic based on a decomposition framework from \cite{luo2022repeated}.
However, all these previous studies in VPP do not require the vehicles to serve all customers in the graph, thus the proposed models and algorithms cannot be applied to the delivery problem directly.

\subsection{Research Gap}

From the reviewed papers, we found that while MV technology has received some attention in public transportation systems, its applications in logistics remain underexplored. Existing studies mainly focus on fixed-route transportation, neglecting the potential of MVs in more flexible, non-fixed-route scenarios. Although a few studies have addressed routing problems involving MVs or vehicle platoons, these are primarily variants of the PDP or the DARP and do not account for VRP settings. Additionally, due to the limited research on the MVRP, there is a lack of discussion on the unique properties of MVs, such as en-route cargo transfer, electricity, and heterogeneous units. Therefore, this study complements existing research by exploring the application of MVs in logistics and addressing these distinctive characteristics.

\section{Problem Definition and Formulation}
\label{sec:problem}

This section defines the problem and constructs the mathematical models. Additionally, we compare the MVRP and VRP solutions theoretically.

\subsection{Problem Definition}

Our problem is defined on a directed graph \( G = (\mathcal{V}, \mathcal{A}) \), where \( \mathcal{V} = \{0, 1, \ldots, N, 0'\} \) represents the set of nodes. $N$ is the maximum number of customers. Node \( 0 \) is the central depot, which serves as the starting node for all feasible routes, while node \( 0' \) is a duplicate of node \( 0 \) and serves as the ending node for all feasible routes. $\mathcal{A}=\{(0, j) \mid j \in N\} \cup\{(i, j) \mid i \neq j\} \cup \{\left(i, 0^{\prime}\right) \mid i \in N \}$ represents the set of arcs. The depot and customers are distributed in a two-dimensional space. The set \( \mathcal{N} = \mathcal{V} \backslash \{0, 0'\} \) denotes the customers, each with a demand \( q_i \). The travel time for each arc \( (i, j) \in \mathcal{A} \) is denoted by \( d_{ij} \), which satisfies the triangle inequality.

\begin{figure*}
    \centering
    \centering
    \includegraphics[width=1.\linewidth]{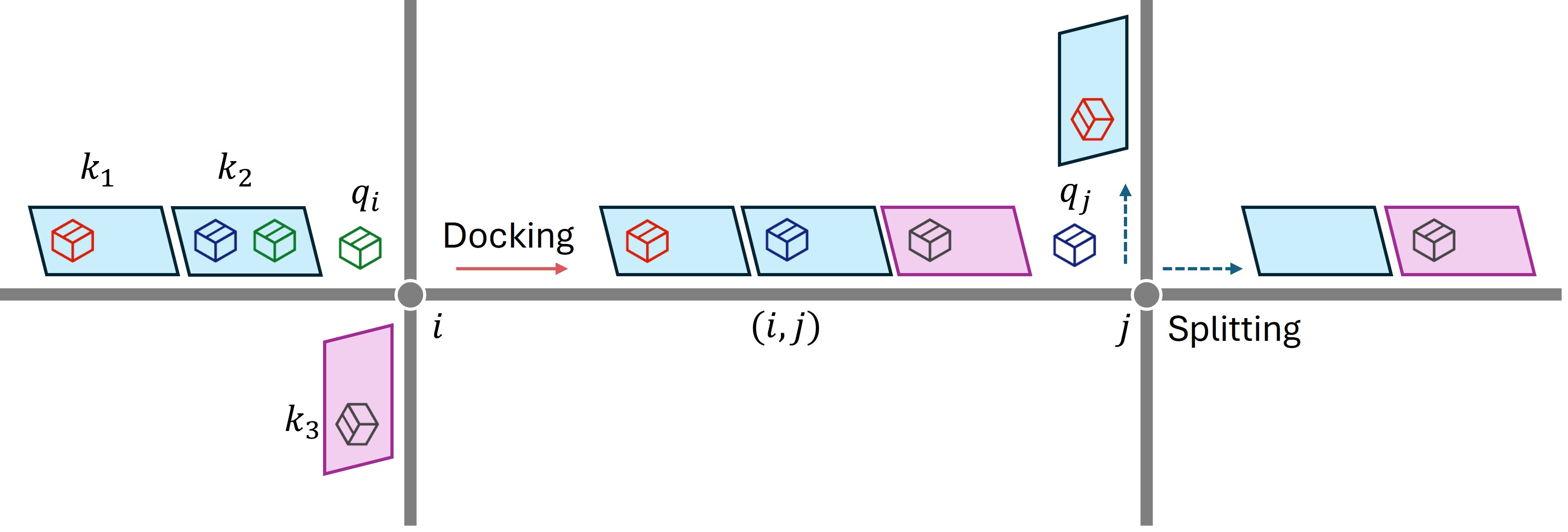}
    \caption{An illustration of the docking and splitting operations for MVs.}
    \label{fig:docking}
\end{figure*}

A fleet of identical $K$ MVs $\mathcal{K}=\{1, \ldots, K\}$ with capacity $Q$ are ready at the depot to serve all customers in $\mathcal{N}$ by platoons. An MV platoon consists of one or more identical connected MV(s). Denote $\mathcal{L}=\{1, 2, ..., L\}$ as the set of the feasible size of an MV platoon, where $L$ is the maximum platoon length. The constraint of the size is named platoon capacity constraint. As illustrated in Figure \ref{fig:docking}, an MV platoon can adjust its size by docking or splitting with other MV(s) at a customer node. In Figure \ref{fig:docking}, two MV platoons meet at customer \(i\) and dock into a single MV platoon consisting of three MVs. After traveling along arc \((i,j)\), this platoon splits at customer \(j\) into two MV platoons. Note that the platoons formed after splitting at \(j\) are different from the two platoons that merged at \(i\). Specifically, in the new platoons, \(k_2\) and \(k_3\) are grouped together, rather than \(k_1\) and \(k_2\). One customer can and can only be served by one MV in the passing platoon, which means that each customer's demand is independent and indivisible during transportation.

\begin{figure*}
   \centering
   \includegraphics[width=.7\textwidth]{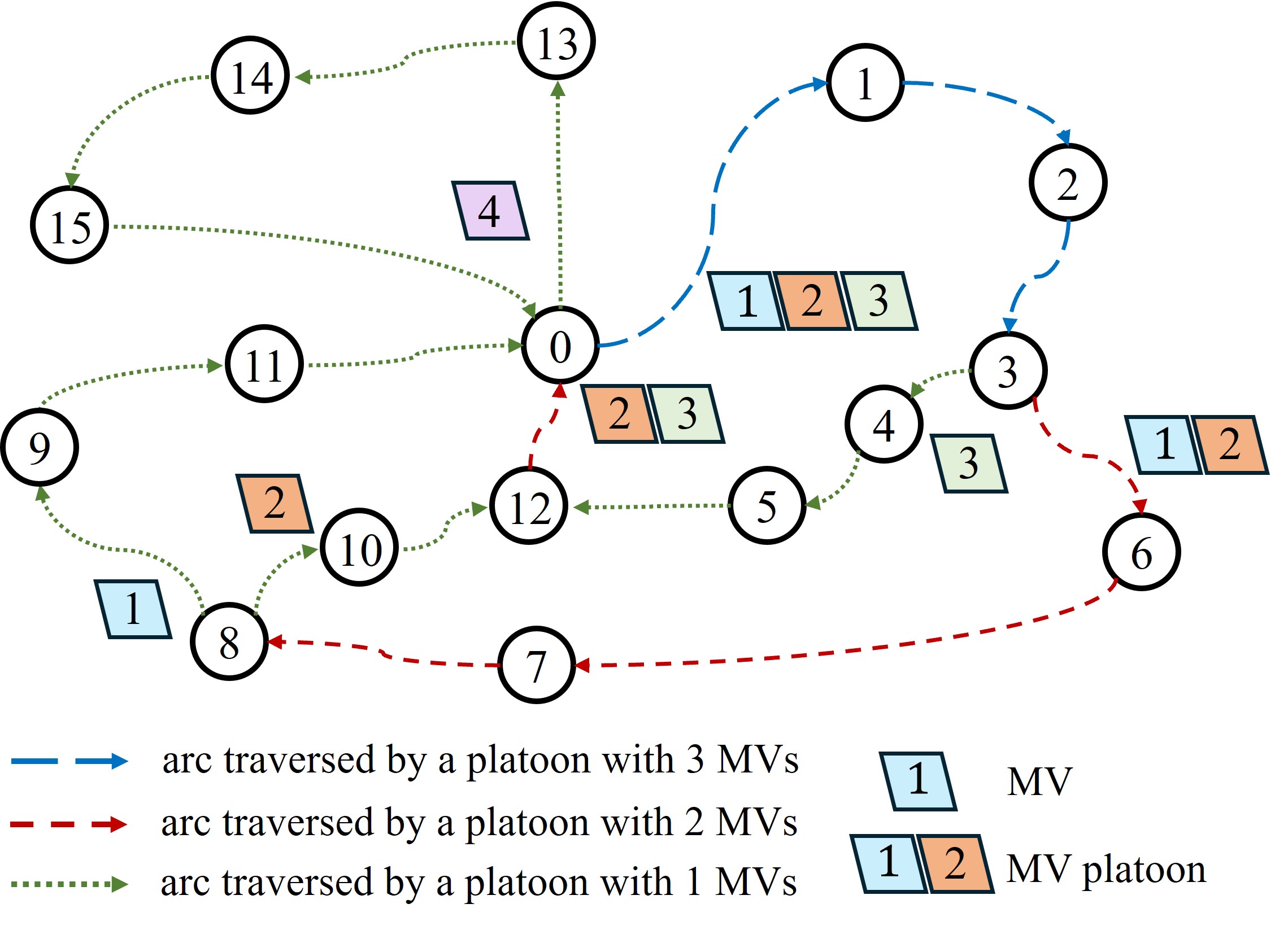}
   \caption{A feasible solution of the MVRP.}
   \label{fig:sol} 
\end{figure*}

Figure \ref{fig:sol} shows a feasible solution of MVRP with 15 customers and 4 MVs. In this solution, a platoon with three MVs, i.e., $k_1$, $k_2$, $k_3$, and $\mathrm{MV}_D$ starts from the depot. This platoon splits at customer 6 into two sub-platoons. After splitting with $k_1$ at customer 8, $k_2$ docks with $k_3$ at customer 12. $\mathrm{MV}_D$ does not dock or split with the other 3 MVs.

In some logistics delivery tasks, customers may have no time window constraints and such problems can be referred to as CMVRP. However, customers often have time window constraints in urban on-demand delivery scenarios, transforming the problem into an MVRPTW. In the MVRPTW, each customer $i$ has a service time $h_i$ and a time window $[e_i, l_i]$, which specifies that the service cannot start either before or after the designated time window. Furthermore, the MVRPTW must consider the time synchronization of MV platoons. Specifically, when multiple MV platoons merge, they must all arrive at the same customer location before the docking process can begin. For example, in Figure \ref{fig:docking}, assuming \(k_1\) and \(k_2\) arrive at customer \(i\) earlier than \(k_3\), they need to wait for \(k_3\) to arrive before performing docking. Waiting is permitted at all locations. Both docking and splitting operations take a fixed reformulation time $t_f$.

The objective of the MVRP is to minimize the total cost of all MVs while satisfying customer demands. Compared to using multiple trucks for transportation, the key advantage of MV platoons lies in their ability to reduce energy consumption per unit distance. Following \cite{fu2023dial}, we adopt a simplified energy cost function with a cost-saving rate $\eta$. The cost for a platoon with $l \in \mathcal{L}$ MVs to traverses arc $(i,j)\in \mathcal{A}$ is defined as $c^l_{ij}=d_{ij} \cdot l \cdot [1-\eta(l-1)]$. This function considers the economies of scale, which means that the longer the platoon, the lower the average operating costs per MV in the platoon. Notice that our formulation and algorithm can be applied to other functions, such as the function proposed in \cite{chen2019operational}.

\subsection{MILP Models}

This section proposes three-index arc-flow MILP models for CMVRP and MVRPTW. We first define the following decision variables:
\begin{itemize}
   \item $x^k_{ij}$: a binary variable equals to 1 if MV $k\in\mathcal{K}$ traverses arc $(i,j)\in \mathcal{A}$;
   \item $y^l_{ij}$: a binary variable equals to 1 if a platoon with $l \in \mathcal{L}$ MVs traverses arc $(i,j)\in \mathcal{A}$;
   \item $u_i$: a decision variable represents the service order of customer $i\in \mathcal{V}$;
   \item $\omega^k_i$: a binary variable equals to 1 if customer $i\in \mathcal{V}$ is served by MV $k\in\mathcal{K}$;
   \item $\delta^k_i$: a decision variable equals to the total demand served by MV $k\in\mathcal{K}$ after arriving customer $i\in\mathcal{V}$.
\end{itemize}

The MILP model for CMVRP is constructed as follows:
\begin{align}
   \min & \sum_{l\in \mathcal{L}}\sum_{(i, j) \in \mathcal{A}}c^l_{ij}y^l_{ij} \label{equ:obj} \\
   \text { s.t. } & \sum_{k\in  \mathcal{K}}\sum_{(i, j) \in \mathcal{A}} x_{i j}^{k}\geq 1 &  \forall j \in \mathcal{N} \label{equ:vis} \\
   & \sum_{(i, j) \in \mathcal{A}} x_{i j}^{k} = \sum_{(j, i) \in \mathcal{A}} x_{j i}^{k} &  \forall k \in \mathcal{K}, j \in \mathcal{N} \label{equ:flow} \\
   & \sum_{(0, j) \in \mathcal{A}} x_{0 j}^{k} =\sum_{(i, 0^{\prime}) \in \mathcal{A}} x_{i 0^{\prime}}^{k}=1 & \forall k \in \mathcal{K} \label{equ:depot} \\
   & u_i + 1 - M^1_i(1-x^k_{ij}) \leq u_j  & \forall (i,j) \in \mathcal{A}, k \in \mathcal{K} \label{equ:st} \\
   & \delta^k_i + q_j - M^2_i(2-x^k_{ij}-\omega^k_j) \leq \delta^k_j  & \forall (i,j) \in \mathcal{A}, k \in \mathcal{K} \label{equ:cap1} \\
   & \delta^k_i - M^3_i(1-x^k_{ij}) \leq \delta^k_j  & \forall (i,j) \in \mathcal{A}, k \in \mathcal{K} \label{equ:cap2} \\
   & \sum_{k\in  \mathcal{K}}\omega^k_i=1 & \forall i \in \mathcal{N} \label{equ:cap3} \\
   & \sum_{(i, j) \in \mathcal{A}} x_{i j}^{k} \geq \omega^k_j & \forall j \in \mathcal{N}, k \in \mathcal{K} \label{equ:cap4} \\
   & \sum_{k\in  \mathcal{K}}x^k_{ij}\leq L &  \forall (i,j) \in \mathcal{A} \label{equ:mpn} \\
   & \sum_{l\in \mathcal{L}} y_{i j}^{l} \geq \sum_{k\in \mathcal{K}}x^k_{ij} / M^4 &  \forall (i,j) \in \mathcal{A} \label{equ:pn} \\
   & y_{i j}^{l} \leq 1 + (l-\sum_{k\in \mathcal{K}}x^k_{ij}) / M^4 &  \forall l\in \mathcal{L}, (i,j) \in \mathcal{A} \label{equ:pn2} \\
   & x_{i j}^{k} \in\{0, 1\} & \forall (i, j) \in \mathcal{A}, k \in \mathcal{K} \label{equ:ijk} \\
   & y_{i j}^{l} \in\{0, 1\} & \forall l\in \mathcal{L}, (i, j) \in \mathcal{A} \label{equ:y} \\
   & \omega_{i}^{k} \in\{0, 1\} & \forall i \in \mathcal{V}, k \in \mathcal{K} \label{equ:omega} \\
   & 0 \leq \delta^k_{i} \leq C &  \forall i \in \mathcal{V}, k \in \mathcal{K} \label{equ:delta} \\
   & u_i \geq 0 & \forall i \in \mathcal{N} \label{equ:u}
\end{align}

The objective function \eqref{equ:obj} minimizes the total cost of all routes. 
Constraints \eqref{equ:vis} ensure that each customer is visited by an MV at least once. 
Constraints \eqref{equ:flow} ensure the flow conservation. 
Constraints \eqref{equ:depot} mean that each MV can leave the depot at most once. 
Constraints \eqref{equ:st} ensure there is no sub-tour in the solution. Here $M^1_i = N$.
Constraints \eqref{equ:cap1}--\eqref{equ:cap3} ensure that vehicle capacity constraints are respected. To be specific, constraints \eqref{equ:cap1} represent that the total demand served by MV $k$ will be increased if the MV traverses from $i$ to $j$ and is responsible for $i$. Here $M^2_i = Q+q_i$. 
Constraints \eqref{equ:cap2} represent that the capacity of MV $k$ will not be reduced if the MV traverses from $i$ to $j$. Here $M^3_i = Q$. 
Constraints \eqref{equ:cap3} ensures that each customer has 1 corresponding MV. 
Constraints \eqref{equ:cap4} restrict that customers can only be assigned to MV that arrive at the customer.
Constraints \eqref{equ:mpn} is the platoon capacity constraint. 
Constraints \eqref{equ:pn}--\eqref{equ:pn2} restrict that $y^l_{ij}$ equals to 1 if the number of MVs in the platoon traverse from $i$ to $j$ is $l$. Here $M^4=L$.
Constraints \eqref{equ:ijk}--\eqref{equ:a} define the domains of decision variables.

For the MVRPTW, we introduce a decision variable to represent the arrival time to replace the decision variable $u_i$:
\begin{itemize}
   \item $a_i^k$: a decision variable which represents the arrival time of MV $k\in \mathcal{K}$ at the node $i\in \mathcal{V}$.
   \item $\alpha_i$: a decision variable which represents the start serving time for the MV platoon at the node $i\in \mathcal{V}$.
\end{itemize}

Constraints \eqref{equ:st} and \eqref{equ:u} are then replaced by the following constraints:
\begin{align}
   & \alpha_i + d_{ij} + h_i - M^5_{ij}(1-x^k_{ij}) \leq a_j^k  & \forall (i,j) \in \mathcal{A}, k \in \mathcal{K} \label{equ:time} \\
   & a_i^k \leq \alpha_i & \forall i \in \mathcal{N}, k \in \mathcal{K} \label{equ:time2} \\
   & e_i \leq a^k_i \leq l_i & \forall i \in \mathcal{N}, k \in \mathcal{K} \label{equ:a} \\ 
   & e_i \leq  \alpha_i \leq l_i & \forall i \in \mathcal{N} \label{equ:alpha}
\end{align}
Constraints \eqref{equ:time} calculate the arrival times of the MV platoon at node $i$, where $M^5_{ij}=\sum_{(i,j) \in \mathcal{A}} d_{ij} + \sum_{i\in \mathcal{N}}s_i$.
Constraints \eqref{equ:time2} represent the synchronization constraints, which require the MV arriving earlier at node $i$ to wait for the later-arriving MV.
Constraints \eqref{equ:a}-\eqref{equ:alpha} are the time window constraints.

\subsection{Theoretical Bound Compared to the VRP}

Finally, we conclude this section by comparing the optimal MVRP solution with the optimal VRP solution. Theorem \ref{the1} offers insights into the differences between these two problems.

\begin{myThe}
   For a given graph and fleet of MVs, denote the optimal VRP solution as $\pi_{VRP}$ and the optimal MVRPTW solution as $\pi_{MVRP}$. $c(\pi)$ is the cost of the solution $\pi$, $c_{max}$ is the maximum average cost for one MV in the unit distance (i.e., the cost with one MV in the platoon), and $c_{min}$ as the minimum average cost for one MV in the unit distance (i.e., the average cost with $L$ MV in the platoon). The following relation holds:
   $$
   \frac{c_{min}}{c_{max}} c(\pi_{VRP}) < c(\pi_{MVRP}) \leq c(\pi_{VRP}).
   $$
   \label{the1}
   \vspace{-0.2cm}
\end{myThe}

\begin{figure*}
   \centering
   \includegraphics[width=1.\textwidth]{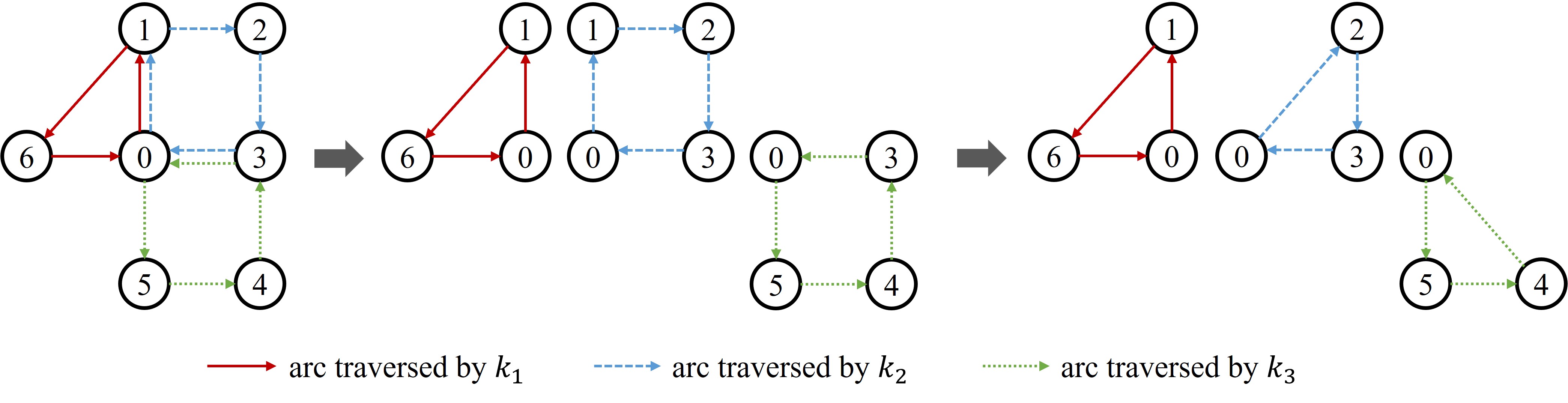}
   \caption{An example of solution separation.}
   \label{fig:the} 
\end{figure*}

\noindent \textit{Proof of Theorem \ref{the1}:} 
The right side of the inequality is easy to understand since all feasible VRP solutions are also feasible for the MVRP. We start from the optimal MVRP solution $\pi_{MVRP}$ to prove the left-side inequality. The solution can be decomposed into a set of routes $\mathcal{R}$ with one MV in each route, as shown in the left side of Figure \ref{fig:the}. For routes $\mathcal{R}$, we have $c(\pi_{MVRP})\geq \sum_{r\in \mathcal{R}} c_{min} d_r$, where $d_r$ is the distance of route $r\in \mathcal{R}$. 

A feasible VRP solution, denoted as $\mathcal{R}'$, can be derived from $\mathcal{R}$ by removing duplicate customers, as illustrated on the right side of Figure \ref{fig:the}. It is worth noting that this step remains valid for MVRPTW because removing customers from a feasible solution does not violate time window constraints. According to the triangle inequality, the cost of the new solution follows $\sum_{r\in \mathcal{R}'} c_{min} d_r \leq \sum_{r\in \mathcal{R}} c_{min} d_r$. Since $\mathcal{R}'$ is not better than the optimal VRP solution, we have $\sum_{r\in \mathcal{R}'} c_{max} d_r \geq s(VRP)$. Combining the three inequalities, we obtain the left-side inequality. \hfill $\square$

Theorem \ref{the1} provides a lower bound of the MVRP solution by solving the VRP. It also reveals that, when the economies of scale disappear, i.e., $c_{max}-c_{min}\to 0$, the problem converges to VRP.

\section{Tabu Search Algorithm}
\label{sec:algorithm}

Although the MILP model can solve the MVRP optimally, this method is very time-consuming. We can easily prove that MVRP is an NP-hard problem since it can be reduced to VRP when the docking and splitting operations are not allowed, or when the maximum number of MVs in a platoon $L$ is set to 1. In this section, we present a TS algorithm to solve large-size MVRP instances in an acceptable time.

\subsection{Solution Representation}

We first introduce the solution representation based on the Gantt chart. In the VRP, a simple array structure is sufficient to describe a single route. However, the MVRP is significantly more complex. Each MV must consider not only its own route and the customers it serves but also its position within an MV platoon. Besides, the time synchronization between MVs is also a significant challenge. These difficulties make it impossible to use a single array to represent the relationships between MVs and MV platoons. To accurately capture the structure of the MVRP, we introduce several new concepts:

\begin{itemize}
   \item \textbf{Route.} The sequence of customers traversed by MV $k\in  \mathcal{K}$.
   \item \textbf{Docking/splitting customers}. Customers at which docking or splitting operations occur.
   \item \textbf{Segment.} The sequence of customers for an MV platoon between two consecutive docking or splitting customers.
   \item \textbf{Path.} The sequence of segments that an MV platoon follows.
\end{itemize}

The foundation of our data structure lies in the concepts of segments and paths. Essentially, a segment and a path correspond to a "node" and a "route" from the perspective of an MV platoon. Using these two concepts, the MVRP solution is divided into two levels: the upper level, which represents the MV platoon perspective, and the lower level, which represents the individual MV perspective. This hierarchical structure can also be leveraged in subsequent algorithms to accelerate the feasibility check and the evaluation of objective values for new solutions.

\begin{figure*}
   \centering
   \includegraphics[width=1.\textwidth]{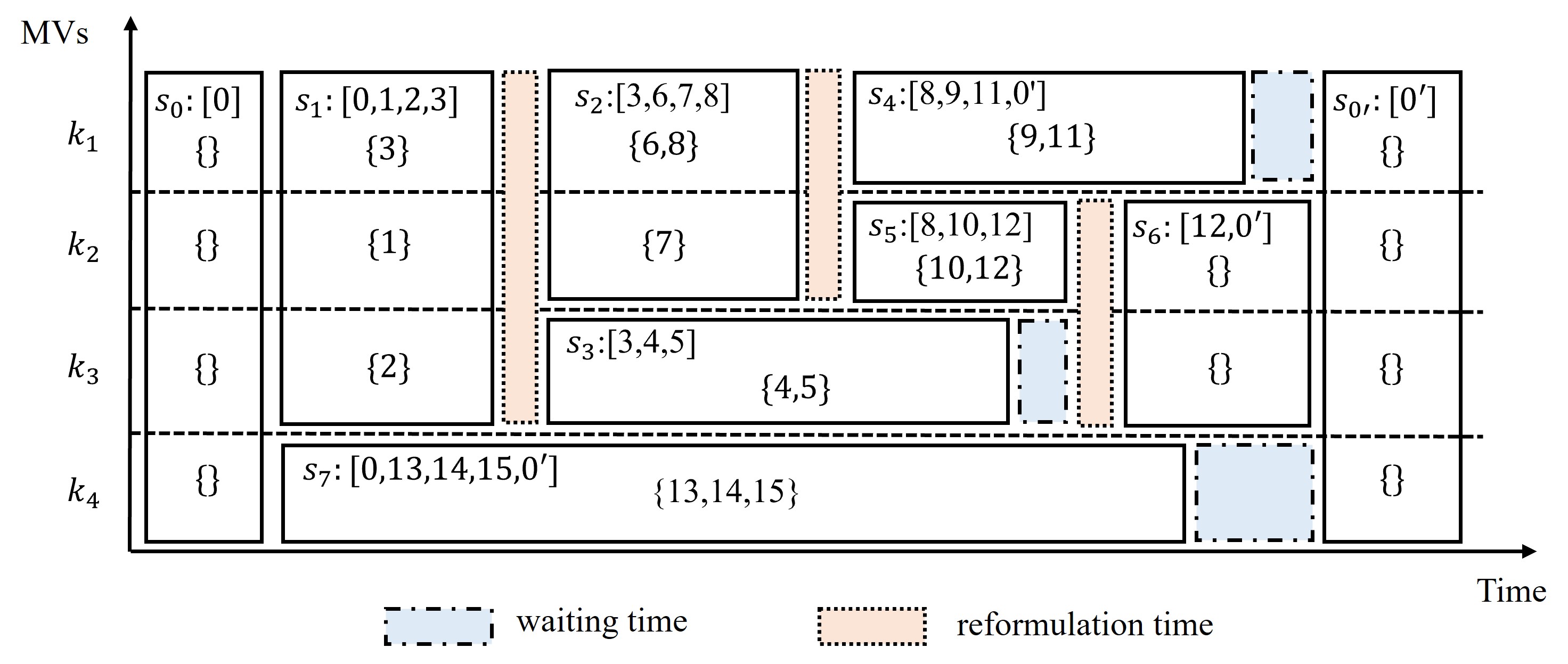}
   \caption{The Gantt chart for solution in Figure \ref{fig:sol}.}
   \label{fig:gantt} 
\end{figure*}

Based on these new structures, the Gantt chart can accurately represent both the cargo flow and the time flow for MV routes and MV platoon paths. The Gantt chart was originally developed as a project management tool to visualize task timelines and was later adapted for applications such as the job shop scheduling problem \citep{jia2007integration}. The Gantt chart corresponding to the solution in Figure \ref{fig:sol} is shown in Figure \ref{fig:gantt} as an example. As illustrated in Figure \ref{fig:gantt}, the horizontal axis represents the arrival time of segments, while the vertical axis represents the sequence of MVs. Denote segment $i$ as $s_i$, where $s_0$ and $s_{0'}$ represent two dummy segments that only contain the depot. Each row corresponds to the route of an MV. For instance, for \(k_1\), the sequence of segments it traverses is \([s_0, s_1, s_2, s_4, s_{0'}]\). Each column in one segment represents the set of customers served by each MV in the platoon. For example, for segment \(s_1\), the associated MVs include \(k_1, k_2\), and \(k_3\), who should serve customer sets $\{3\}, \{1\}$, and $\{2\}$, respectively. Segments are required to record the sequence of visited customers as well as the preceding and succeeding segments. By maintaining information about adjacent segments, the time synchronization between MVs can be evaluated by comparing the arrival time and departure time of all succeeding and preceding customers, as illustrated by the waiting time and reformulation time in Figure \ref{fig:sol}. In the CMVRP, where time flow is not considered, the width of each segment can be set to a unit length.

\subsection{Neighborhood Structure}

TS is a memory-based heuristic search strategy that guides local search beyond local optimality \citep{glover1989tabu, glover1990tabu}. In the TS algorithm, the most critical component is the neighborhood structure, i.e., the operators. Common operators used in VRP include relocation, 2-opt, and cross-exchange \citep{mcnabb2015testing}. However, existing operators designed for VRP cannot be directly applied to MVRP due to differences in solution structures. Moreover, these VRP-based neighborhoods do not encompass MVRP solutions, making it impossible to search for MVRP solutions effectively. In this section, we design two merging operators that can transform VRP solutions into MVRP solutions and adapt the famous relocation operator for the MVRP.

\subsubsection{Merging Operators.}

\begin{figure*}
   \centering
   \includegraphics[width=.9\textwidth]{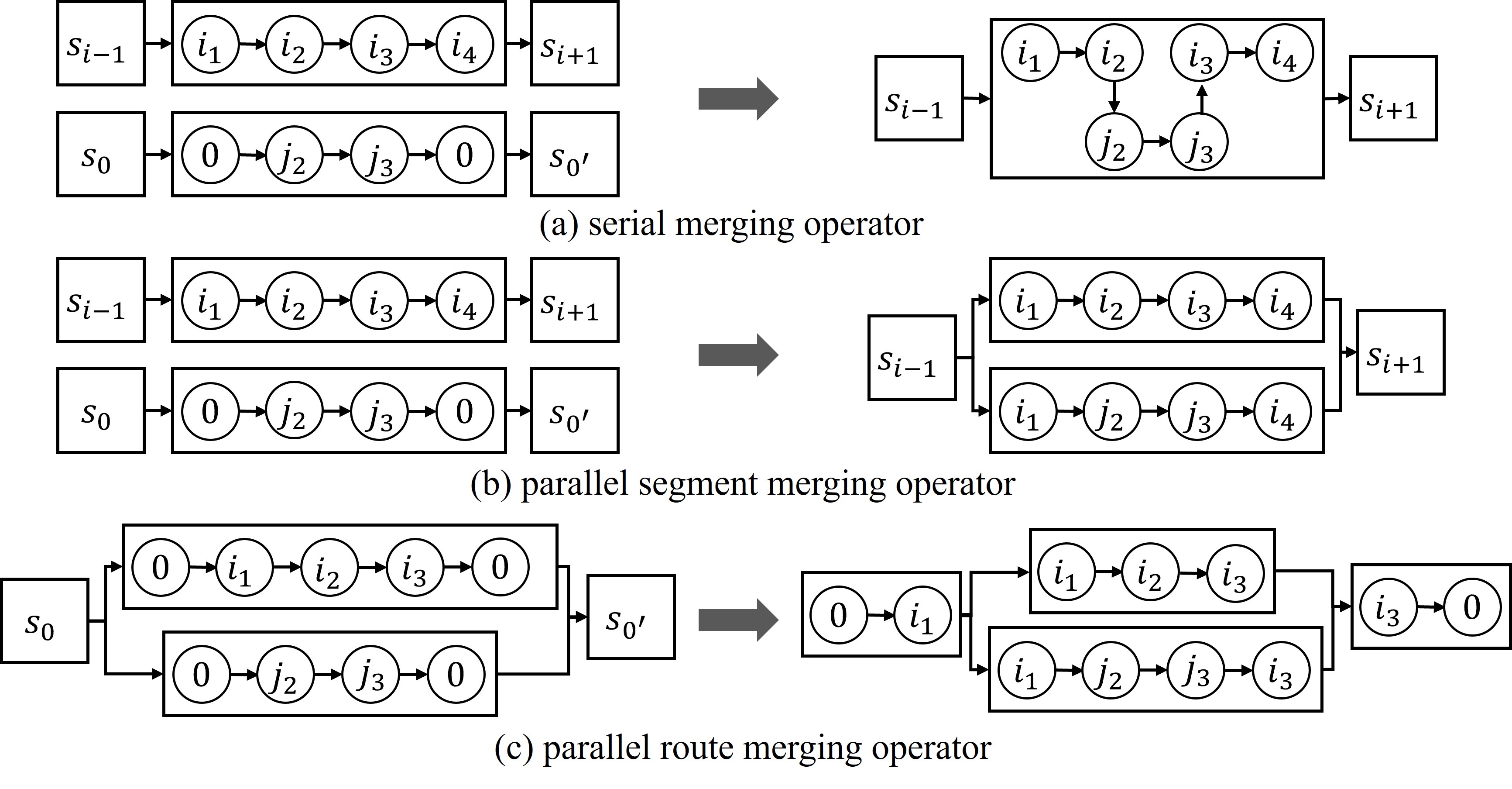}
   \caption{Examples of the merging operators.}
   \label{fig:cmo} 
\end{figure*}

This section introduces three types of merging operators tailored for the MVRP. Figure \ref{fig:cmo} illustrates three examples of these operators, where rectangles represent segments and circles represent customers. Merging operators are typically applied between segments. For simplicity, we use the term \textit{depot segment} to denote a segment that starts and ends at the depot in the following section.

The first merging operator is the serial merging operator (SMO). The SMO inserts all customers from a depot segment $s_j$ into a specific position $p$ in another segment $s_i$. As illustrated in Figure \ref{fig:cmo}(a), all customers in segment \( s_j \) are inserted into segment \( s_i \) after $i_2$ to form a new segment. In the TS algorithm, the neighborhood of SMO includes all segments $s_i$, depot segments $s_j$, and insertion positions $p$.

Different from the sequential merging in SMO, the parallel segment merging operator (PRMO) attempts to integrate the depot segment $s_j$ into a segment $s_i$ to form an MV platoon. As illustrated in Figure \ref{fig:cmo}(b), $s_i$ and $s_j$ are merged into a single MV platoon, where the first node \(i_1\) in segment \(s_i\) is designated as the splitting customer, and the last node \(i_S\) is designated as the docking customer. PMO also allows $s_j$ to merge with $s_i$ either with the path before $s_i$ or in the path after $s_i$. Therefore, the neighborhood of PMO includes all segments $s_i$, depot segments $s_j$, and two binary values representing the merging options for the forward path and the afterward path.

The third type of merging operator is the parallel route merging operator (PRMO), which is designed for two depot segments, $s_i$ and $s_j$. Similar to the PMO, these two depot segments can also be merged into an MV platoon. However, this operator offers greater flexibility in selecting the positions of the splitting and docking customers. Figure \ref{fig:cmo}(c) illustrates an example where \(i_1\) and \(i_3\) are selected as the splitting and docking customers, respectively. The neighborhood of this operator includes all pairs of depot segments, $s_i$ and $s_j$, along with the positions of the splitting and docking customers in segment $s_i$.

In the SMO and PSMO, When the segments' platoons comprise multiple MVs, the algorithm should combine as many MVs as possible to reduce the number of MVs in the new platoon, which reduces the cost. An enumeration method with a time complexity of $O(n!)$ can be used for the pairing procedure. To do it more efficiently, we propose a pairing algorithm with a time complexity of $O(n)$, which is described in Algorithm \ref{alg:pair}. This algorithm satisfies Corollary \ref{cor:pair}.

\begin{algorithm}[htbp]
   \small
      \renewcommand{\algorithmicrequire}{\textbf{Input:}}
      \renewcommand{\algorithmicensure}{\textbf{Output:}}
      \caption{The pseudo-code of the pairing algorithm.}
      \label{alg:pair}
      \begin{algorithmic}[1]
         \STATE \textbf{Input:} MV set $\mathcal{K}_{j}$ for depot segment $s_j$ and $\mathcal{K}_{i}$ for segment $s_i$.
            \STATE Sort $\mathcal{K}_{j}$ and $\mathcal{K}_{i}$ in non-ascending and non-descending order concerning the demand of each MV, respectively;
            \STATE Initialize an empty set of MV pairs $\mathcal{P}$;
            \FORALL{MV $k_1 \in \mathcal{K}_{j}$}
            \STATE Denote $k_2$ as the first MV in $\mathcal{K}_{i}$;
            \IF{$\delta_{k_1}+\delta_{k_2} \leq Q$}
            \STATE $\mathcal{P} \leftarrow \mathcal{P} \cup \{(k_1, k_2)\}$, $\mathcal{K}_{i} \leftarrow \mathcal{K}_{i} \backslash \{k_2\}$;
            \ENDIF
            \ENDFOR
            \STATE \textbf{Output:} $\mathcal{P}$.
      \end{algorithmic}  
\end{algorithm}
   
\begin{myCor}
   Algorithm \ref{alg:pair} can always obtain the optimal pairs (i.e., minimize the number of MVs after pairing).
   \label{cor:pair}
\end{myCor}

\noindent \textit{Proof of Corollary \ref{cor:pair}:}
We prove that Algorithm \ref{alg:pair} can always obtain the optimal pair solution by induction. Denote $f(\mathcal{K}_j, \mathcal{K}_i)$ as the maximum number of pairs between $\mathcal{K}_j$ and $\mathcal{K}_i$. Let $k_1$ represent the MV in $\mathcal{K}_j$ with the maximum demand and $k'_1$ represent the MV in $\mathcal{K}_i$. First, we prove that if $k_1$ cannot pair with $k'_1$, then dropping $k_1$ will not alter the optimal number of pairs. Since $k'_1$ has the minimum demand, $k_1$ cannot pair with any MV in $\mathcal{K}_i$. Therefore, $f(\mathcal{K}_j \backslash {k_1}, \mathcal{K}_i) = f(\mathcal{K}_j, \mathcal{K}_i)$. Next, we prove that if $k_1$ can pair with $k'_1$, the number of pairs remains unchanged if we pair $k_1$ with $k'_1$. Assume the optimal pair solution contains $(k_1, k'_2)$ and $(k_2, k'_1)$, where $k_2$ and $k'_2$ are two MVs in $\mathcal{K}_j$ and $\mathcal{K}_i$, respectively. The pairs $(k_1,k'_1)$ and $(k_2,k'_2)$ are also feasible. Thus, we have $f(\mathcal{K}_j \backslash {k_1}, \mathcal{K}_i \backslash {k'_1})+1 = f(\mathcal{K}_j, \mathcal{K}_i)$. Therefore, each iteration in steps 4-9 will not affect the optimal number of pairs. \hfill $\square$

When there are MVs in the segments that cannot be merged, it is necessary to determine how each MV will travel from depot 0 to $s_i$ and from $s_i$ back to depot 0'. To minimize the cost of adding an MV to that segment, the algorithm needs to search for the shortest paths without violating the platoon capacity constraint. Based on the Gantt chart data structure, a breadth-first search algorithm is implemented in our TS algorithm.

\subsubsection{Relocation Operators.}

\begin{figure*}
   \centering
   \includegraphics[width=0.9\textwidth]{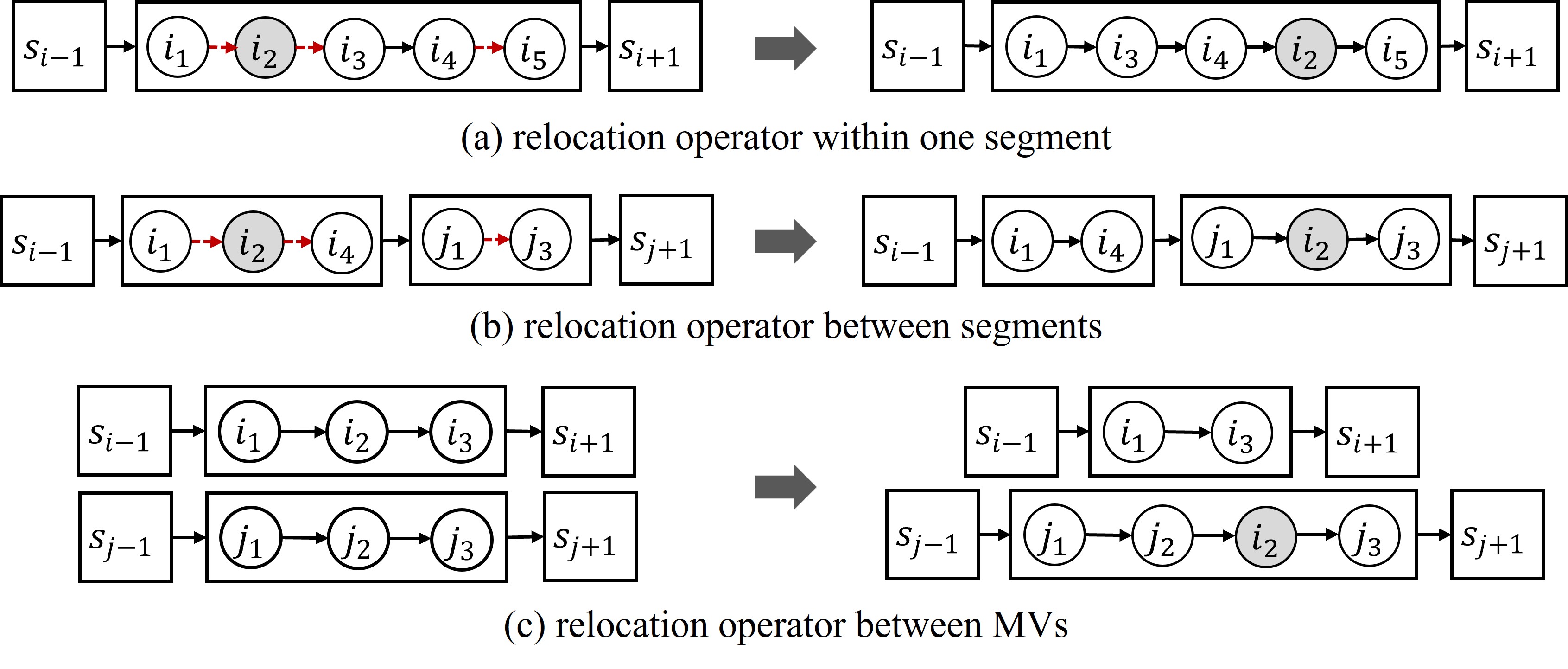}
   \caption{Examples of the relocation operators.}
   \label{fig:re} 
\end{figure*}

The relocation operators relocate a customer from the original position to a new position. According to the solution structure of the MVRP, we introduce three relocation operators, illustrated in Figure \ref{fig:re}:
\begin{itemize}
   \item Relocate a customer into a new position in the same segment and the same MV.
   \item Relocate a customer to a new position in another segment, where the original and new segments are all served by the same MV.
   \item Relocate the customer from an MV to another MV.
\end{itemize}

In the TS algorithm, a tabu list is used to prevent redundant solution searches. The tabu list implemented in our algorithm is a two-dimensional array corresponding to each arc. Each element in the array indicates the number of iterations during which the corresponding arc is prohibited. For the relocation operators, the red dashed lines in Figures \ref{fig:re} represent the arcs that are marked as tabu for each operator. Since the merging operators do not result in redundant solutions, the tabu list is not applied to them.

\subsubsection{Acceleration for Time Window Constraints.}

Although both merging operators and relocation operators can be applied to CMVRP and MVRPTW, many solutions in the neighborhood of MVRPTW may violate time window constraints. Checking time window constraints requires calculating the arrival time for each customer in the new solution. This process is computationally expensive in MVRPTW because MVs need to remain synchronized during docking and splitting. As a result, even a minor modification may affect the arrival times of customers in other MV platoons.

To accelerate the time window constraint verification, we propose a method based on calculating the latest arrival time. Consider a segment $s=[i_0,i_1,...,i_I]$, where $I$ is the maximum number of customers in $s$. For a given customer $i$ with arrival time $a_i$, we can compute the maximum delay at this node as $\tau_i = l_i - a_i$. This implies that the MV can be delayed by at most $\tau_i$ from the depot while still satisfying the time window constraint at customer $i$. Therefore, by calculating backward, we can compute the latest arrival time at customer $j$ as
\begin{align}
    \eta_j = a_j - \max_{j\leq i\leq I} \{\tau_i\}
\end{align}
to account for all subsequent time window constraints in the segment. When an operator generates a new solution, it is only necessary to compute the arrival time of the first unchanged customer and compare it with its latest arrival time to verify whether the time window constraint is satisfied. This method is also applicable to VRPTW.

\begin{figure}
    \centering
    \includegraphics[width=0.9\linewidth]{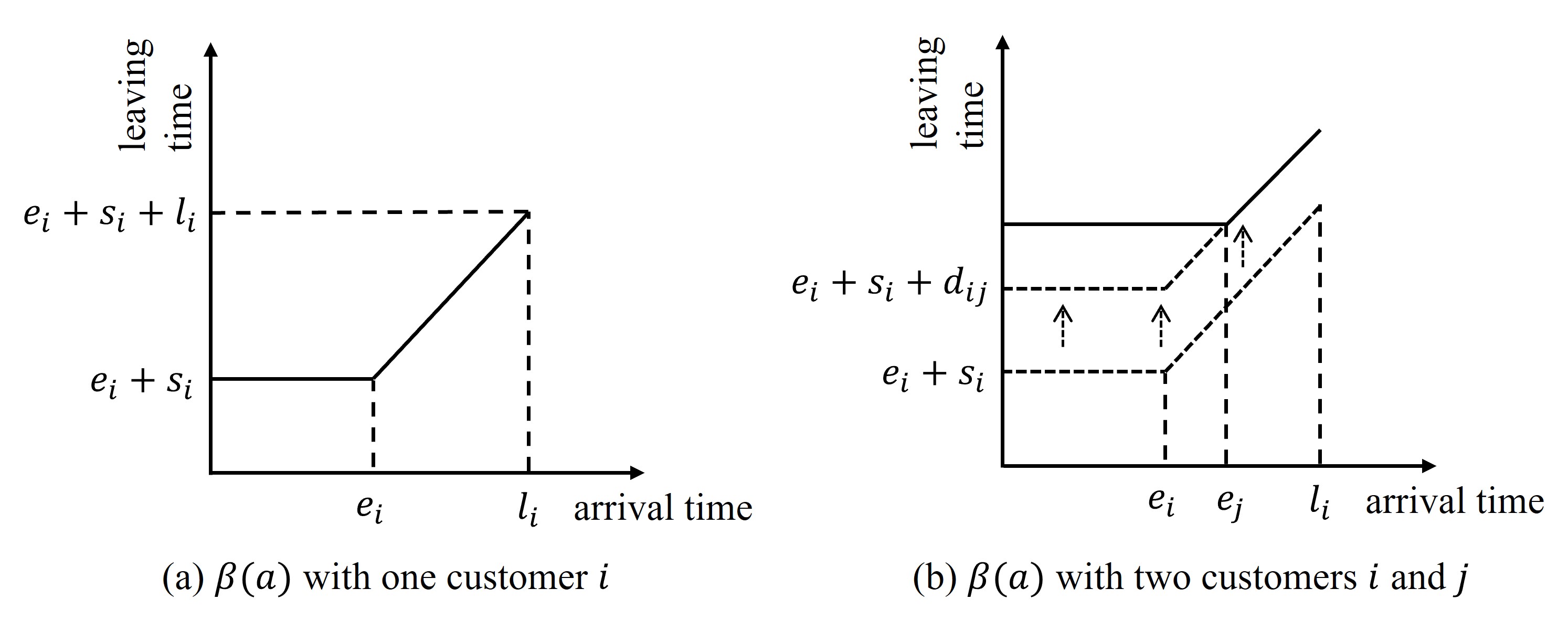}
    \caption{Examples of function forms of $\beta(a)$ with one and two customers.}
    \label{fig:beta}
\end{figure}

For segments with unchanged visiting sequences, this method can also be used to compute the latest arrival time of the entire segment. However, when considering a path that includes multiple segments, the problem becomes more complex. Consider a path $[s_0,s_1,s_2,s_3,s_{0'}]$. Suppose an operator modifies $s_1$, causing the arrival time at the end of $s_1$ to increase by $\Delta a$. Due to MV platoon synchronization and waiting for time windows to open, the increased arrival time at $s_3$ may not equal $\Delta a$.

To address this issue, we define a function $\beta(a)$ to compute the departure time from a segment as a function of its arrival time $a$. We first consider a case with a single customer $i$. Under time window constraints, the departure time from $i$ can be computed as $\beta_i = \max\{e_i, a\} + s_i$, where the function form is illustrated in Figure \ref{fig:beta}(a). For two customers, the departure time from $j$ can be calculated as $\beta_j = \max\{e_j, \beta_i + d_{ij}\} + s_j$. The function form is shown in Figure \ref{fig:beta}(b). We observe that the function remains piecewise, with a constant segment followed by an affine segment with slope 1. It is straightforward to prove that this result holds for more customers. Therefore, we only need to compute the coordinates of the breakpoints to determine the functional form of $\beta$. Since the affine function maintains identical $x$- and $y$-axis coordinates, we only need to compute the departure time corresponding to $a=0$.

\subsection{Shaking}

If the best solution has not been improved for a certain number of iterations, shaking is applied to the current solution to escape from the local optimum \citep{brandao2011tabu}. The idea of shaking is to separate several MV routes in the solution as depot segments. An advantage of this shaking method is that it remains valid for MVRPTW because removing customers from a feasible solution does not violate time window constraints. Combined with the merging operators, this step can effectively expand the search space for MVRP solutions.

\subsection{Generate Initial Solutions}

As a local search algorithm, TS often gets trapped in a local optimum. Therefore, a multi-start strategy is applied to the TS \citep{marti2013multi}. The idea is to sample a solution pool \(\mathcal{S}\) as the initial set of solutions for the TS. After the optimization of TS, the best solution is selected as the final output.

The Clarke–Wright savings algorithm (CW algorithm) proposed in \cite{clarke1964scheduling} is applied to obtain the initial solutions. The algorithm initially considers each customer as a separate route directly connected to the depot. It then calculates the 'savings' that would be achieved by merging any two of these individual routes, where 'savings' refers to the reduced cost. The merging operator with the highest savings is applied to the solution. This method iteratively reduces the total number of routes and the overall distance until no further savings can be achieved.

In our algorithm, we first randomly sample a certain number of sparse VRP solutions. 'Sparse' means that the number of customers in each route is significantly less than the capacity constraint limit. The CW algorithm then uses the merging operators to optimize these solutions. Finally, these solutions are fed into the TS for further optimization.

\section{Discussion on Variants of the MVRP}
\label{sec:variants}

Although this paper focuses on the fundamental CMVRP and MVRPTW, the MVRP can also be extended to various other practical application scenarios. This section explores several extensions of the MVRP based on the potential capabilities and application scenarios of MVs mentioned in the literature. Specifically, we discuss three types of variants: MVRP with en-route cargo transfer, electric MVRP, and MVRP with Heterogeneous modular units.

\textbf{MVRP with en-route cargo transfer.} Allowing the en-route cargo transfer between MVs in the same MV platoon can provide more flexibility for the routing problem. Similar to modular buses, where passengers can move between MV units \citep{tian2022planning}, MV can be designed with automated sliding doors at the junctions between adjacent MVs. By using automated conveyor belts or manually transferring by onboard couriers, cargo can be transferred between MVs in the same MV platoon.

\begin{figure*}
   \centering
   \includegraphics[width=.5\textwidth]{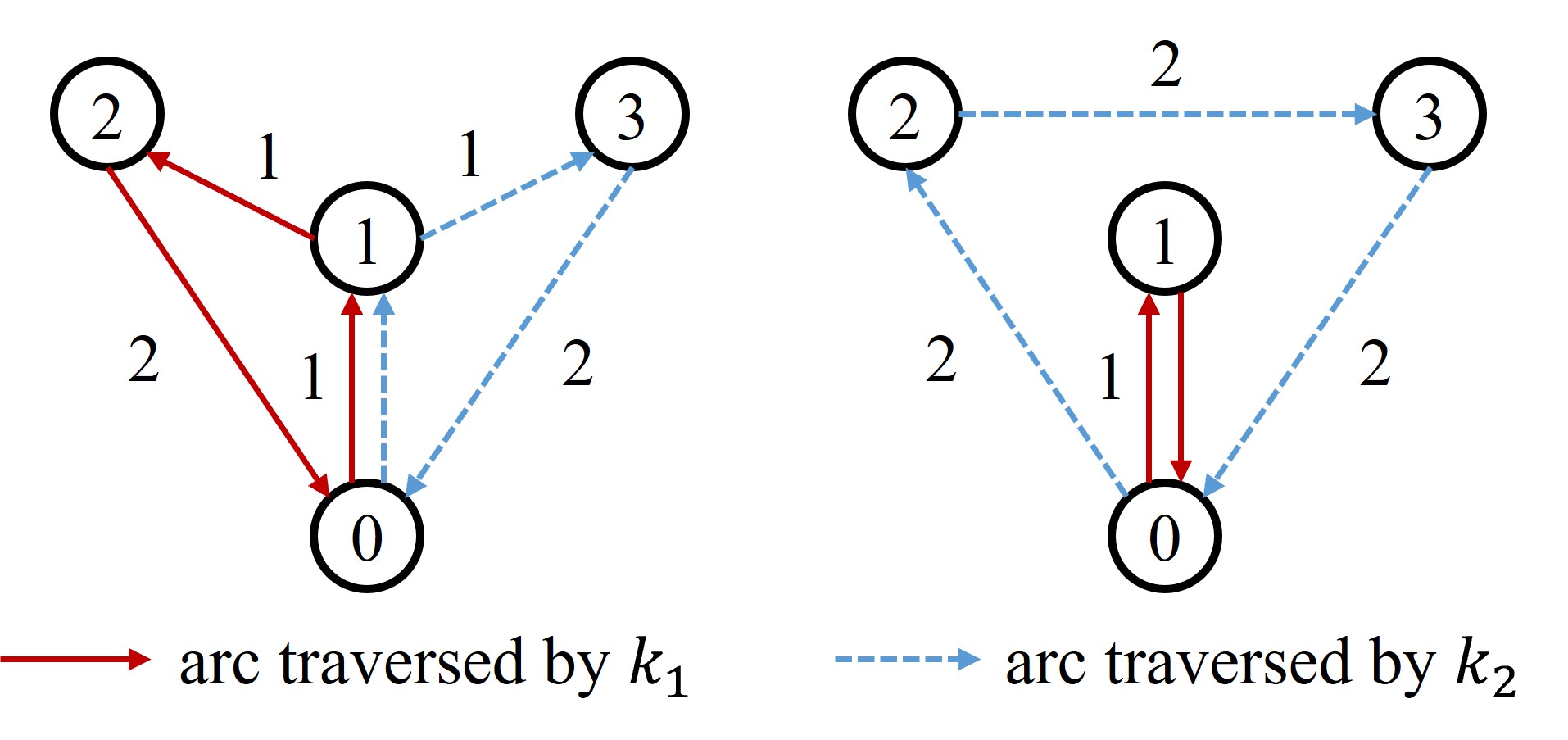}
   \caption{Example of en-route cargo transfer.}
   \label{fig:transfer} 
\end{figure*}

The benefit of the en-route transfer is shown in Figure \ref{fig:transfer}. In this graph, we set $d_1=2,d_2=1,d_3=1,Q=2$, and $c^1_{ij}=d_{ij},c^2_{ij}=1.5d_{ij}$. The red solid line and blue dashed line are the routes of the two MVs, respectively. If the en-route transfer is not allowed, the only optimal VRP and MVRP solution will be the same, as shown in the right figure with a cost equal to 8. However, if en-route transfer is allowed, one can serve customer 1 and receive 1 demand transfer from another MV, as shown in the left figure. This solution has a lower cost 7.5 than the optimal solution of the MVRP without en-route cargo transfer.

\textbf{Electric MVRP.} 
Electric vehicles (EVs) have the potential to reduce the environmental impact of transportation, which has been widely discussed in the literature \citep{kucukoglu2021electric}. Combining the technologies of EVs and MVs, the Electric MVRP extends the routing problem by incorporating energy constraints and charging requirements. In the Electric MVRP, MVs can be designed to support energy sharing when connected as an MV platoon, potentially balancing battery levels across units to extend their operational range. Furthermore, individual MV units can be dispatched to charging stations while the remaining units continue deliveries or prepare for subsequent tasks. After recharging, the recharged unit can rejoin the platoon, ensuring minimal disruption to the overall schedule.

\textbf{MVRP with heterogeneous modular units.} 
MVs also have the capability to transport different types of demands. MV platoons can be formed using multiple types of modules, such as freight and passenger units, to achieve integration of diverse transportation tasks \citep{hatzenbuhler2023modular,ulrich2019new}. This flexibility is particularly valuable in logistics operations that require specialized handling for certain freight. For example, perishable items, such as fresh food, often necessitate cold chain logistics with temperature-controlled modules to ensure quality and safety during transportation \citep{deng2023pickup,qin2019vehicle}. The MVRP with heterogeneous modular units incorporates different types of modules and can therefore meet diverse demands within the same transportation network.

\section{Computational Experiments}
\label{sec:experiments}

This section presents the instance sets and analyzes the performance of the algorithm. Our algorithm was coded in Java programming language using ILOG CPLEX 12.6.3 as the solver. The experiments were conducted on a machine equipped with a 3.2GHz AMD Ryzen 7 7735HS with Radeon Graphics CPU and 16G of memory under the Windows 11 operating system.

\subsection{Instances}

The benchmark instances used in our experiments are slightly modified from well-known CVRP and VRPTW benchmarks in the literature. The CMVRP instances are derived from Augerat Set A and Set B \citep{augerat1995approche}. To demonstrate the impact of instance size on the problem, we randomly selected instances from these sets and removed some nodes. As a result, we obtained ten instances with 10 nodes as small-size instances and five instances each with 30, 50, and 67 nodes as large-size instances (since we could not find at least 5 instances with more than 67 nodes). The MVRPTW instances are based on the Solomon benchmark instances \citep{solomon1987algorithms}. Similar to the CMVRP instances, we selected 5 instances with 10 nodes as small-size instances and 5 instances with 50 nodes as large-size instances. The instances are labeled as ``setID\_$N$\_ID", where the set ID represents A and B for Augerat Set A and Set B in CMVRP, and C for the Solomon benchmark in MVRPTW. The number of nodes and the instance ID are also included in the labels. For small-size instances, the vehicles' capacity is modified to 70, and the maximum platoon length \( L \) is set to 2. For large-size instances, the capacity is 100, and \( L \) is 3. Additionally, to simulate road conditions in urban delivery environments, we use the Manhattan distance to calculate the distance between nodes. Similar to \cite{fu2023dial}, we set \(\eta=0.1\) in the cost function. The code for all algorithms and data analysis, as well as the benchmark instances, are publicly available at \url{url} (will update after the acceptance of this paper).

\subsection{Algorithm Performance}

In this section, we first evaluate the algorithm's performance on small-size instances by comparing it with the MILP model from Section \ref{subsub:cplex}. Next, we assess its performance on large-size instances in Section \ref{subsub:algorithm}. Finally, we evaluate the algorithm's components, including neighborhood operators and strategies, in Section \ref{subsub:components}.

\subsubsection{Performance Comparison between the TS Algorithm and the MILP.}
\label{subsub:cplex}

\begin{table}[htbp]
  \centering
  \caption{Comparision on MILP model and the heuristic algorithm on small-size instances.}
    \begin{tabular}{lrrrrrrrr}
    \toprule
    \multirow{2}[4]{*}{Instance} & \multicolumn{5}{c}{MILP}              & \multicolumn{2}{c}{TS} & \multirow{2}[4]{*}{$\Delta$} \\
    \cmidrule{2-8}          & $UB$    & $LB$    & $Gap$ & $Time$ & $\#Node$  & $Obj$    & $Time$  &  \\
    \midrule
    A-10-1 & 541.8 & 438.6 & 19.04 & 1800  & 1753653 & 541.8 & \textbf{0.40} & 0.00 \\
    A-10-2 & 451.0 & 451.0 & 0.00  & 282   & 152317 & 451.0 & \textbf{0.15} & 0.00 \\
    A-10-3 & 447.2 & 447.2 & 0.00  & 97    & 39894 & 447.2 & \textbf{0.12} & 0.00 \\
    A-10-4 & 357.2 & 198.6 & 44.39 & 1800  & 917360 & \textbf{356.2} & \textbf{0.16} & -0.28 \\
    A-10-5 & 420.0 & 420.0 & 0.00  & 193   & 25327 & 420.0 & \textbf{0.11} & 0.00 \\
    B-10-1 & 427.8 & 234.5 & 45.18 & 1800  & 1801329 & 427.8 & \textbf{0.05} & 0.00 \\
    B-10-2 & 336.0 & 336.0 & 0.00  & 55    & 17696 & 336.0 & \textbf{0.46} & 0.00 \\
    B-10-3 & 510.8 & 510.8 & 0.00  & 970   & 594673 & 510.8 & \textbf{0.22} & 0.00 \\
    B-10-4 & 510.8 & 469.5 & 8.09  & 1800  & 511880 & 510.8 & \textbf{0.12} & 0.00 \\
    B-10-5 & 452.0 & 279.7 & 38.13 & 1800  & 459082 & \textbf{448.0} & \textbf{0.23} & -0.89 \\
    C-10-1 & \textbf{354.8} & 267.2 & 24.68 & 1800  & 877931 & 355.4 & \textbf{0.26} & 0.17 \\
    C-10-2 & 89.2  & 89.2  & 0.00  & 362   & 153625 & 89.2  & \textbf{0.16} & 0.00 \\
    C-10-3 & 227.0 & 227.0 & 0.00  & 3     & 8074  & 227.0 & \textbf{0.21} & 0.00 \\
    C-10-4 & 347.4 & 164.1 & 52.75 & 1800  & 1914134 & \textbf{346.0} & \textbf{0.19} & -0.40 \\
    C-10-5 & 367.4 & 367.4 & 0.00  & 51    & 20605 & 367.4 & \textbf{0.05} & 0.00 \\
    \midrule
    Average & 389.4 & 326.7 & 15.48 & 974   & 616505 & 389.0 & 0.19  & -0.09 \\
    \bottomrule
    \end{tabular}%
  \label{tab:milp}%
\end{table}%

Table \ref{tab:milp} reports the results of the MILP model and heuristic algorithm on small-size instances, where the time limit on each run of the MILP model was set to 1800 seconds. In the following sections, we use columns $UB$ and $LB$ to represent the best upper and lower bounds for the MILP model. Column $\#Node$ is the number of nodes explored in the B\&B (branch and bound) tree in the MILP model. Column $Gap$ is the percentage difference between the best upper bounds and the lower bound. Column $Time$ reports the time in seconds consumed to solve one instance. Specifically, $Gap=(UB-LB)/UB \times 100$. Column $Obj$ represents the objective values of the TS algorithm. Column $\Delta$ is the percentage gap of the objective values between the MILP model and the TS algorithm, calculated as $\Delta=(Obj-UB)/Obj \times 100$. 

From Table \ref{tab:milp} we can notice that the MILP model can only solve 8 out of 15 instances with 10 nodes to optimality within the time limit. However, the average computation time for the 10 instances is 974 seconds. Moreover, the optimality gaps are, on average, 15.48\%, which is still large. In contrast, the TS algorithm significantly reduces the computation time while maintaining solution quality. The TS algorithm found all the solutions obtained by the MILP model and two better solutions for instances "A-10-5", "B-10-5", and "C-10-4", where the MILP did not find the optimal solutions. All instances were solved within 0.50 seconds, with an average solution time of 0.19 seconds. This demonstrates that our TS algorithm outperforms the MILP model.

\subsubsection{Performance of the TS Algorithm on Large-size Instances.} 
\label{subsub:algorithm}

\begin{table}[htbp]
  \centering
  \caption{Performance of the TS algorithm on large-size instances.}
    \begin{tabular}{lrrlrrlrrlrr}
    \toprule
    Instance & $Obj$  & $Time$  & Instance & $Obj$  & $Time$  & Instance & $Obj$   & $Time$ & Instance & $Obj$   & $Time$ \\
    \midrule
    A-30-1 & 784.0 & 5.1   & A-50-1 & 1255.4 & 8.6   & C-50-1 & 1001.2 & 12.3  & A-67-1 & 1458.8 & 13.2 \\
    A-30-2 & 843.2 & 2.5   & A-50-2 & 1374.8 & 8.5   & C-50-2 & 782.2 & 7.2   & A-67-2 & 1981.0 & 13.7 \\
    A-30-3 & 869.6 & 2.5   & A-50-3 & 1250.6 & 7.8   & C-50-3 & 978.0 & 10.0  & B-67-1 & 1294.6 & 12.4 \\
    B-30-1 & 667.8 & 3.1   & B-50-1 & 950.8 & 12.4  & C-50-4 & 925.4 & 17.2 & B-67-2 & 1482.0 & 13.4 \\
    B-30-2 & 798.2 & 3.4   & B-50-2 & 1633.6 & 7.9   & C-50-5 & 942.4 & 10.6  & B-67-3 & 1360.0 & 12.5 \\
    \midrule
    Average & 792.6 & 3.3   & Average & 1293.0 & 9.0   & Average & 925.8 & 11.5  & Average & 1515.3 & 13.0 \\
    \bottomrule
    \end{tabular}%
  \label{tab:heuristic}%
\end{table}%

\begin{figure}
    \centering
    \includegraphics[width=0.8\linewidth]{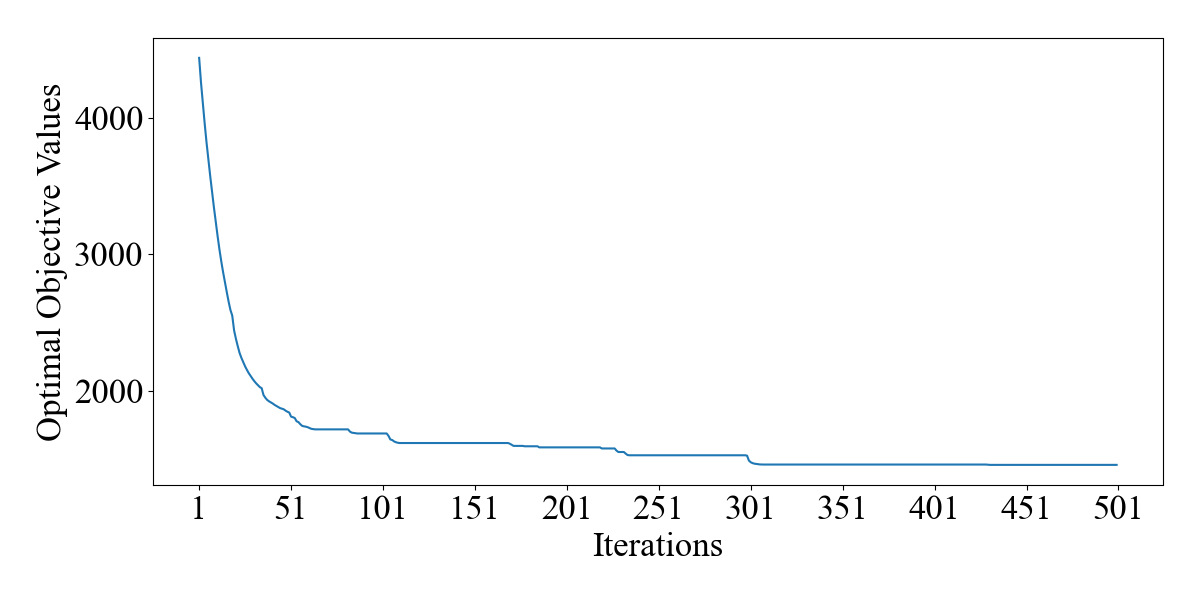}
    \caption{Convergence plots for the TS algorithm for instances “A-67-1”.}
    \label{fig:con}
\end{figure}

Table \ref{tab:heuristic} shows that our algorithm can solve large-size instances at a fast computation speed. For instances with 30, 50, and 67 nodes, the algorithm completes within approximately 5, 10, and 14 seconds, respectively. We further provide the convergence plot for instance “A-67-1” in Figure \ref{fig:con}, which illustrates the algorithm's rapid convergence in the first few iterations. The plot indicates the fast convergence of the algorithm in the first 50 iterations. After about 300 iterations, the solution converges to the optimal value.

\subsubsection{Performance of Components in the Modified TS.}
\label{subsub:components}

\begin{table}[htbp]
  \centering
  \caption{Comparision on algorithms without each component.}
    \begin{tabular}{lrrrrr}
    \toprule
    Instance & $Obj$  & $Obj_{R}$ & $Obj_{M}$ & $Obj_{S}$ & $Obj_{MS}$ \\
    \midrule
    A-30-1 & 784.0 & 828.4 & 828.4 & 828.4 & 828.4 \\
    A-30-2 & 843.2 & 890.8 & 878.8 & 878.8 & 878.8 \\
    A-30-3 & 869.6 & 923.4 & 890.0 & 909.4 & 888.0 \\
    B-30-1 & 667.8 & 746.6 & 744.6 & 744.6 & 700.0 \\
    B-30-2 & 798.2 & 812.6 & 804.6 & 804.6 & 804.6 \\
    \midrule
    Average & 792.6 & 840.4 & 829.3 & 833.2 & 820.0 \\
    \midrule
    A-50-1 & 1255.4 & 1297.6 & 1297.6 & 1297.6 & 1297.6 \\
    A-50-2 & 1374.8 & 1398.8 & 1378.8 & 1378.8 & 1374.8 \\
    A-50-3 & 1250.6 & 1280.4 & 1272.4 & 1272.4 & 1250.6 \\
    B-50-1 & 950.8 & 959.2 & 957.2 & 957.2 & 951.6 \\
    B-50-2 & 1633.6 & 1704.4 & 1678.4 & 1678.4 & 1678.4 \\
    \midrule
    Average & 1293.0 & 1328.1 & 1316.9 & 1316.9 & 1310.6 \\
    \midrule
    C-50-1 & 1001.2 & 1089.2 & 1024.0 & 1083.2 & 1083.2 \\
    C-50-2 & 782.2 & 792.2 & 782.2 & 786.2 & 786.2 \\
    C-50-3 & 978.0 & 1089.2 & 1024.0 & 1083.2 & 1024.0 \\
    C-50-4 & 925.4 & 942.4 & 932.0 & 940.4 & 940.4 \\
    C-50-5 & 942.4 & 976.0 & 940.0 & 960.4 & 960.4 \\
    \midrule
    Average & 925.8 & 977.8 & 940.4 & 970.7 & 958.8 \\
    \midrule
    A-67-1 & 1458.8 & 1475.4 & 1471.4 & 1471.4 & 1463.4 \\
    A-67-2 & 1981.0 & 2034.4 & 2013.2 & 2013.2 & 1981.0 \\
    B-67-1 & 1294.6 & 1308.0 & 1300.0 & 1300.0 & 1300.0 \\
    B-67-2 & 1482.0 & 1506.0 & 1498.0 & 1498.0 & 1482.0 \\
    B-67-3 & 1360.0 & 1410.0 & 1360.0 & 1360.0 & 1360.0 \\
    \midrule
    Average & 1515.3 & 1546.8 & 1528.5 & 1528.5 & 1517.3 \\
    \bottomrule
    \end{tabular}%
  \label{tab:ope}%
\end{table}%

To evaluate the neighborhood operators, i.e., the relocate operators and merge operators, as well as the shaking and multi-start strategies proposed in our algorithm, we compare the objective values obtained by the algorithm without each of these components on large-size instances. Columns $Obj$, $Obj_{R}$, $Obj_{M}$, $Obj_{S}$, and $Obj_{MS}$ in Table \ref{tab:ope} represent the objective values with all components included and the results without the relocation operators, merge operators, shaking, and multi-start strategy. Notably, the algorithm that includes all components performs the best across all instances. 

\begin{figure}
    \centering
    \includegraphics[width=0.8\linewidth]{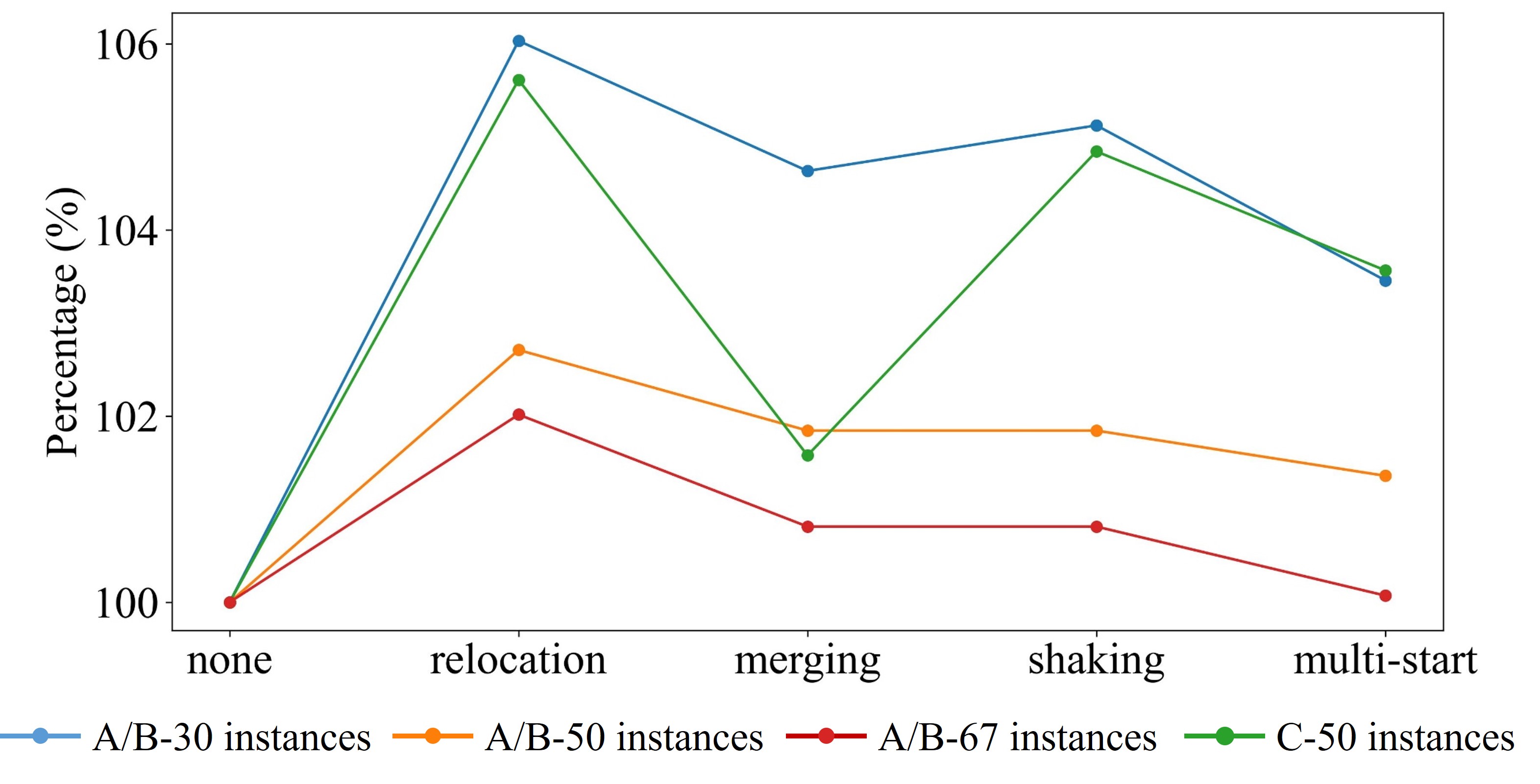}
    \caption{The ratio of results from algorithms without each component relative to the algorithm with all components.}
    \label{fig:com}
\end{figure}

Figure \ref{fig:com} demonstrates the impact of each component with different instance sizes. It can be observed that removing any single component degrades the solution quality produced by the algorithm. This result highlights that each part of the tailored algorithm effectively contributes to improving solution quality. Moreover, among these components, removing the relocation operator and the shaking strategy has the most significant impact on solution quality. This finding further emphasizes the importance of a well-designed neighborhood structure and strategies for escaping local optima in TS. Finally, the figure also illustrates that for small-size instances, each component significantly improves solution quality. However, as the instance size increases, the improvement from each strategy diminishes, indicating that the algorithm’s solution quality is less robust for larger instances.

\subsection{Value of MV Platooning}


\begin{figure}
    \centering
    \includegraphics[width=0.8\linewidth]{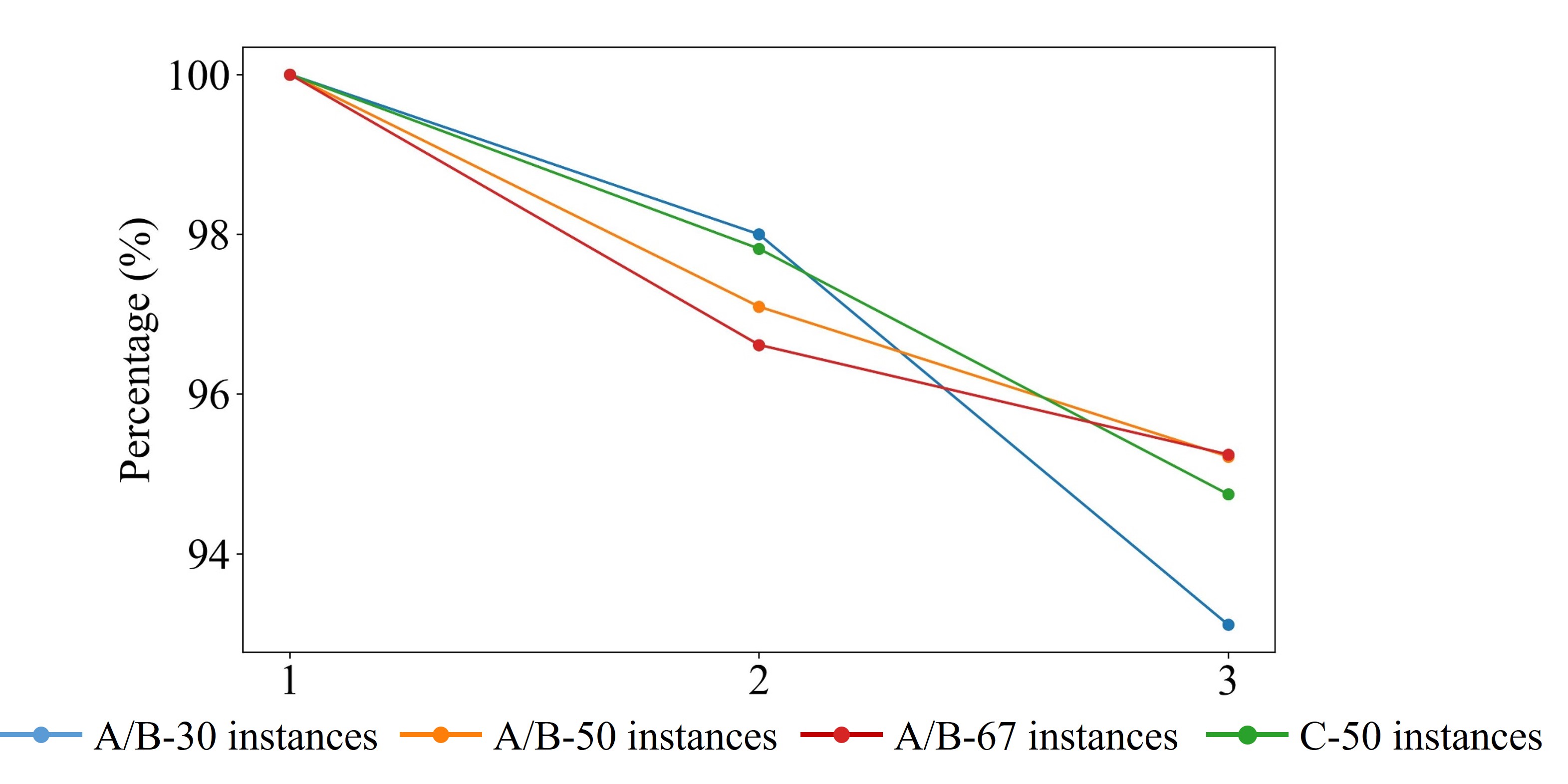}
    \caption{The impact of maximum platoon length $L$.}
    \label{fig:len}
\end{figure}

To quantify the benefits of MV platooning, we conduct experiments with large-size instances to compare the objective values of the MVRP under different maximum platoon lengths, \(L = 1, 2,\) and 3. Figure \ref{fig:len} illustrates the percentage ratio of average objective values for all instances with \(L=1, 2,\) and 3, relative to the result with \(L=1\). The results show that as \(L\) increases, the objective values for all instances decrease. Specifically, when \(L\) increases from 1 to 2, the objective value decreases by approximately 2.6\%. Similarly, when \(L\) increases from 1 to 3, the objective value decreases by about 5.6\%. According to Theorem \ref{the1}, with \(\eta = 0.1\), the theoretical upper bounds for the reduction in objective value are 10\% and 20\% for \(L = 2\) and \(L = 3\), respectively. These results demonstrate that the actual reduction in objective values is approximately one-fourth of the theoretical upper bound.

\subsection{Sensitivity Analysis}

\begin{figure}
    \centering
    \includegraphics[width=0.8\linewidth]{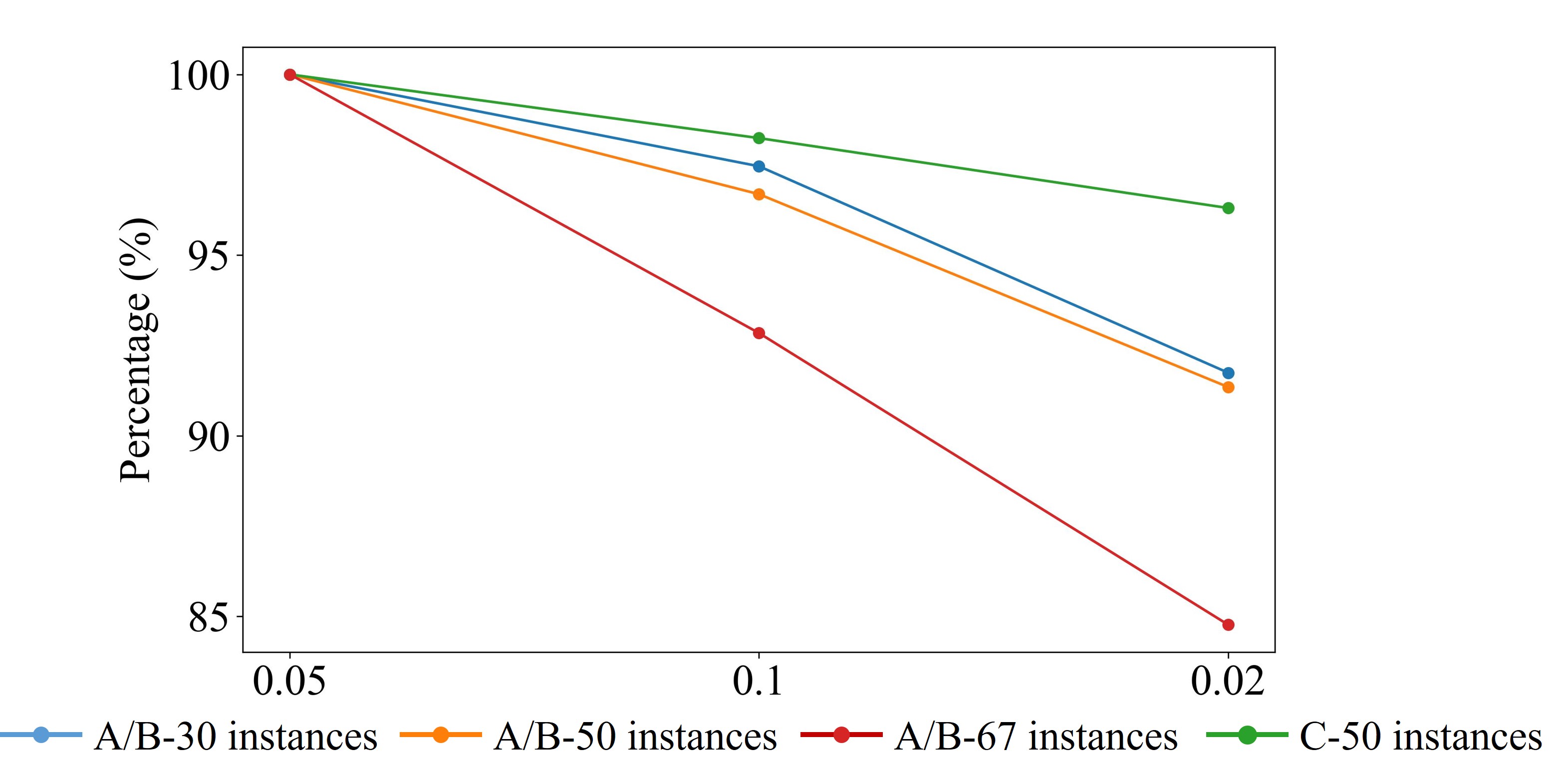}
    \caption{The impact of cost-saving rate $\eta$.}
    \label{fig:eta}
\end{figure}

Since the MVRP is a new variant of the VRP, we conduct a sensitivity analysis for the key parameter $\eta$. As analyzed in Theorem \ref{the1}, the major advantage of the MVRP over the VRP lies in the reduced travel cost when multiple MVs dock to form a platoon. If $\eta = 0$, the optimal solution of the MVRP will be no better than that of the VRP. To investigate the influence of this parameter, we analyze its impact by varying its value from 0.1 to 0.05 and 0.2 in large-size instances. Figure \ref{fig:eta} plots the average results compared to the result with $\eta = 0.05$. We observe that a higher rate leads to lower average costs. Specifically, when $\eta$ increases from 0.05 to 0.1, the cost is reduced by approximately 4\%. Further increasing $\eta$ from 0.05 to 0.2 results in an approximate 8\% cost reduction. These results suggest that improving the cost-saving capability of MV technology through platooning can enhance overall benefits and promote the broader adoption of MVs.

\section{Conclusion}
\label{sec:conclusion}

This paper introduces a new variant of the VRP that incorporates MVs' unique ability to dock and split in a trip. Unlike previous studies, we explore the application of MVs in logistics applications. An MILP formulation is first developed to model the problem, which can be solved using advanced solvers. To efficiently solve large-size instances, we designed a multi-start TS algorithm. A specialized solution representation based on the ST-DAG and Gantt chart is developed to capture the MVRP solution, and multiple local search operators are tailored specifically for the MVRP. To avoid local optima, our TS begins with an initial solution pool generated by the CW algorithm and perturbs the current solution if no improvement is observed after a set number of iterations. 

The results demonstrate that our TS algorithm can find all optimal solutions obtained by the MILP model in small-size instances and exhibit good convergence speed in large-size instances. Furthermore, MV platooning can significantly reduce the total cost by about 5.6\%, aligning with theoretical analysis. However, the effectiveness of our algorithm on larger instances remains uncertain. Evaluating its performance on large-scale cases will require the development of faster exact algorithms beyond the current MILP model.

In addition to the basic MVRP, this study also explores several variants of the MVRP. While our heuristic algorithm can be applied to these variants, the increased problem complexity may require further refinement of the algorithmic details. Future research could focus on developing more efficient solution methods for these variants, such as ALNS, Genetic Algorithms, and Deep Reinforcement Learning. Moreover, exploring the properties of these problems through numerical experiments and theoretical analysis could provide valuable insights. Additionally, future research could investigate efficient exact algorithms for solving the MVRP, including branch-and-price, branch-and-cut, and Benders decomposition techniques.

\bibliographystyle{informs2014}
\bibliography{reference.bib}

\end{document}